\documentclass[12pt]{article}

\usepackage{amsfonts,graphicx,euscript}

\newtheorem{lemma}{Lemma}

\newcommand{\ip}[2]{\langle#1,#2\rangle}
\newcommand{\homo}[1]{\mathbf{#1}}

\title{On semidefinite representations of plane quartics}

\author{Didier Henrion$^{1,2}$}

\begin{document}

\maketitle

\footnotetext[1]{LAAS-CNRS, University of Toulouse, France}
\footnotetext[2]{Faculty of Electrical Engineering,
Czech Technical University in Prague, Czech Republic}
\addtocounter{footnote}{2}

\begin{abstract}
This note focuses on the problem of representing convex
sets as projections of the cone of positive semidefinite
matrices, in the particular case of sets generated by
bivariate polynomials of degree four.
Conditions are given for the convex hull of a plane quartic
to be exactly semidefinite representable with at most 12 lifting
variables. If the quartic is rationally parametrizable,
an exact semidefinite representation with 2 lifting variables
can be obtained.
Various numerical examples illustrate
the techniques and suggest further research directions.
\end{abstract}

\begin{center}\small
{\bf Keywords}\\[1em]
semidefinite programming; polynomials; algebraic plane curves
\end{center}

\section{Introduction}

Following the pioneering work of Nesterov and Nemirovskii
\cite{nn}, it is known that a significant collection of convex sets are
representable as linear sections or projections of
the cone of positive semidefinite matrices (or semidefinite
cone for short). A classification has been proposed in \cite{bn},
and advanced optimization modeling software such as YALMIP \cite{yalmip}
exploits this knowledge to generate semidefinite programming problems
from general convex programming problems, thus allowing the use
of powerful interior-point algorithms based on self-concordant
barrier functions \cite{nn}. The dictionnary proposed in \cite{bn}
is far from being comprehensive, however, and further efforts are
required to identify convex sets which can be efficiently
represented via the semidefinite cone. This note focuses on
the particular case of the convex hull of algebraic plane curves
of degree four.

Let $x \in {\mathbb R}^2 \mapsto p(x)$ be a bivariate
polynomial of degree 4, a quartic, and let
\begin{equation}\label{curve}
{\mathcal C} = \{x \in {\mathbb R}^2 :\: p(x) = 0\}
\end{equation}
be the corresponding algebraic curve.
Alternatively, we can use
the homogenization $\homo{x} \in {\mathbb P}^2 \mapsto \homo{p}(\homo{x})$,
defined in the projective plane ${\mathbb P}^2$, and such that
$\homo{p}(\homo{x})
= \homo{x}_0^{\deg p} p(\homo{x}_1/\homo{x}_0,
\homo{x}_2/\homo{x}_0)$ and $p(x)
= \homo{p}(1,x_1,x_2)$.
Algebraic curve (\ref{set}) can be defined equivalently in the projective plane as
\[
{\mathcal P} = \{\homo{x} \in {\mathbb P}^2 :\:
\homo{p}(\homo{x}) = 0\}.
\]
See the appendix for more information on algebraic sets in the projective plane.

Let
\begin{equation}\label{set}
{\mathcal P} = \mathrm{conv}\:{\mathcal C}
\end{equation}
denote the convex hull of curve $\mathcal C$,
a semi-algebraic set. We assume that $\mathcal P$
has a non-empty interior, and we denote by $\partial\mathcal P$
the boundary of $\mathcal P$. For notational convenience, we also
assume that the sign of $p(x)$ is such that the semialgebraic set
$\{\homo{x} \in {\mathbb P}^2 \: :\: \homo{p}(\homo{x}) \geq 0\}$
is included in $\mathcal P$.

As stated above, we are interested in modeling $\mathcal P$
as a projection of a linear section of the semidefinite cone:
\[
\{\homo{x} \in {\mathbb P}^2 :\:
\exists\: \homo{y} \in {\mathbb P}^m :\: \sum_{i=0}^2 F_i \homo{x}_i
+ \sum_{j=0}^m G_j \homo{y}_j \succeq 0\}
\]
where $F_i, G_j$ are real symmetric
matrices of size $n$ to be found, and $\succeq 0$ means
positive semidefinite. Following the
terminology of \cite{bn}, we are seeking a
semidefinite representation of $\mathcal P$.
The $\homo{y}_j$ are called lifting variables, or liftings
for short, and
the set is obtained by projecting the feasible
set of a linear matrix inequality (LMI).

Note that lifting variables are necessary in the above
description, already in the particular case when
$\{\homo{x} \in {\mathcal P}^2 :\:
\homo{p}(\homo{x}) \geq 0\}$ is a convex set.
Indeed, a generic line passing through set $\mathcal P$ cuts
its boundary $\mathcal C$ at real points only twice, which is
less than the degree of $\homo{p}$, and this implies that
$\mathcal P$ is not rigidly convex in the
sense of \cite{hv}. It follows that $\mathcal P$
cannot be represented without
lifting variables ($m=0$), or, equivalently, as a linear
section of the semidefinite cone.

\section{Primal moment approximations}\label{section_primal}

In \cite{l}, Jean-Bernard Lasserre proposed a hierarchy of semidefinite,
or LMI relaxations
for polynomial optimization, based on the theory of
moments and the dual representation of non-negative polynomials
as sum-of-squares (SOS). In particular \cite{l06},
Lasserre's relaxation of order $k$ for set $\mathcal P$
is given by
\[
{\mathcal P}_k = \{\homo{x} \in {\mathbb P}^2 \: :\: \exists\:
\homo{y} \in {\mathbb P}^m \: :\:
\homo{y}_{00} = \homo{x}_0, \: \homo{y}_{10} = \homo{x}_1, \: \homo{y}_{01} = \homo{x}_2, \: M_k(\homo{y}) \succeq 0,
\: M_{k-2}(\homo{p}\homo{y}) = 0\}
\]
with
\[
\homo{y} = [\:\homo{y}_{00} \:\: \homo{y}_{10} \:\: \homo{y}_{01} \:\: \homo{y}_{20} \:\: \homo{y}_{11}
\:\: \homo{y}_{02} \:\: \homo{y}_{30} \:\cdots\: \homo{y}_{0,2k}]
\in {\mathbb P}^m
\]
a (truncated) vector of bivariate
moments with $m+1=(k+1)(2k+1)$ entries,
\[
M_k(\homo{y}) =  \left[\begin{array}{cccccc}
\homo{y}_{00} & \homo{y}_{10} & \homo{y}_{01} & \homo{y}_{20} \\
\homo{y}_{10} & \homo{y}_{20} & \homo{y}_{11} & \homo{y}_{30} \\
\homo{y}_{01} & \homo{y}_{11} & \homo{y}_{02} & \homo{y}_{21} \\
\homo{y}_{20} & \homo{y}_{30} & \homo{y}_{21} & \homo{y}_{40} \\
& & & & \ddots \\
& & & & & \homo{y}_{0,2k}
\end{array}\right]
\]
a (truncated) moment matrix, and $M_k(\homo{p}\homo{y})$
a (truncated) localising matrix.
For example, if $p(x) = 1-2x^3_1x_2$
we have $\homo{p}(\homo{x})
= \homo{x}^4_0-2\homo{x}^3_1\homo{x}_2$ and
\[
M_k(\homo{p}\homo{y}) = \left[\begin{array}{cccccc}
\homo{y}_{00}-2\homo{y}_{31} & \homo{y}_{10}-2\homo{y}_{41} & \homo{y}_{01}-2\homo{y}_{32} \\
\homo{y}_{10}-2\homo{y}_{41} & \homo{y}_{20}-2\homo{y}_{51} & \homo{y}_{11}-2\homo{y}_{42}  \\
\homo{y}_{01}-2\homo{y}_{32} & \homo{y}_{11}-2\homo{y}_{42} & \homo{y}_{02}-2\homo{y}_{33} \\
& & & \ddots \\
& & & & \homo{y}_{0,2k}-2\homo{y}_{3,1+2k}
\end{array}\right].
\]
In particular if $k=2$ (the first relaxation in the hierarchy for
quartic sets), the localising matrix
is the scalar linear form $M_0(\homo{p}\homo{y})
= \homo{y}_{00}-2\homo{y}_{31}$ obtained by lifting
monomials of $p$.

Note also that by symmetry there are redundant constraints in
the equation $M_{k-2}(\homo{p}\homo{y})=0$ arising in the definition
of ${\mathcal P}_k$. For notational convenience, we however
stick with this matrix notation.

By definition, ${\mathcal P} \subset \cdots \subset {\mathcal P}_3
\subset {\mathcal P}_2$, and hence we have a hierarchy of
embedded outer
semidefinite representable approximations for $\mathcal P$.

\section{Dual SOS approximations}

Let
\[
{\mathcal F} = \{\homo{f} \in {\mathbb P}^2 \: :\:
\homo{f}(\homo{x}) = \homo{f}_0 \homo{x}_0 + \homo{f}_1 \homo{x}_1 +
\homo{f}_2 \homo{x}_2 \geq 0 \:\: \forall \homo{x} \in {\mathcal P}\}
\]
denote the dual cone of $\mathcal P$
in the sense that if $\homo{f} \in {\mathcal F}$ then
the half-plane $\{\homo{x} \: :\: \homo{f}(\homo{x}) \geq 0\}$ contains $\mathcal P$.
It follows that either $\homo{x} \in {\mathcal P}$
or there exists $\homo{f} \in {\mathcal F}$ such that $\homo{f}(\homo{x}) < 0$.
Since $\mathcal P$ has a non-empty interior,
these are mutually exclusive statements.
A geometric interpretation is that $\mathcal P$
is the intersection of the (generally infinite number of)
half-planes generated by all separating elements, i.e.
${\mathcal P} = \{\homo{x} \in {\mathbb P}^2 \: :\:
\homo{f}(\homo{x}) \geq 0 \:\: \forall \homo{f} \in {\mathcal F}\}$.
Set $\mathcal F$ admits the following equivalent representation,
dehomogenized with respect to $x$:
\[
{\mathcal F} = \{\homo{f} \in {\mathbb P}^2 \: :\:
f(x) = \homo{f}_0 + \homo{f}_1 x_1 +
\homo{f}_2 x_2 \geq 0 \:\: \forall x \in {\mathcal P}\}.
\]

Let ${\mathbb R}[x]_k$ denote the set of polynomials
of degree at most $k$ and ${\mathbb S}[x]_{2k}$ the set of polynomials
that can be written as sums of squares (SOS) of polynomials
of ${\mathbb R}[x]_k$. Let
\begin{equation}\label{fk}
{\mathcal F}_k = \{\homo{f} \in {\mathbb P}^2 \: :\:
f(x) = s_0(x) + s_1(x)p(x),
\: s_0(x) \in {\mathbb S}[x]_{2k}, \: s_1(x)
\in {\mathbb R}[x]_{2(k-2)}\}.
\end{equation}
Since SOS polynomials are only a subset of non-negative polynomials
\cite{r}, it holds
${\mathcal F} \supset \cdots \supset
{\mathcal F}_3 \supset {\mathcal F}_2$.

\begin{lemma}\label{lemma_dual}
${\mathcal P}_k = \{x \in {\mathbb R}^2 \: :\: f(x) \geq 0 \:\:
\forall f \in {\mathcal F}_k\}.$
\end{lemma}

{\bf Proof:}
The condition that an element $x^* \in {\mathbb R}_2$, or, equivalently,
an element $\homo{x}^* \in {\mathbb P}_2$, belongs
to ${\mathcal P}_k$ is the existence of a vector $\homo{y}$ satisfying
the primal semidefinite constraints
\begin{equation}\label{primal}
\begin{array}{l}
M_k(\homo{y}) = \sum_{\alpha} A_{\alpha} \homo{y}{\alpha}  \succeq 0 \\
M_{k-2}(\homo{p}\homo{y}) = \sum_{\alpha} B_{\alpha} \homo{y}_{\alpha} = 0 \\
\homo{y}_{00} = \homo{x}^*_0, \: \homo{y}_{10} = \homo{x}^*_1, \: \homo{y}_{01} = \homo{x}^*_2. \\
\end{array}
\end{equation}
Given two symmetric matrices $A$ and $X$ of the same size,
define the inner product $\ip{A}{X} = \mathrm{trace}(AX)$.
Build the Lagrangian $L(\homo{y},X,\homo{f})=
-\ip{M_k(\homo{y})}{X_0}-\ip{M_{k-2}(\homo{p}\homo{y})}{X_1}
+\homo{f}_0(\homo{y}_{00}-\homo{x}^*_0)+\homo{f}_1(\homo{y}_{10}-\homo{x}^*_1)+
\homo{f}_2(\homo{y}_{01}-\homo{x}^*_2)$
with $X_0 \succeq 0$ and the corresponding
dual function $\inf_{\homo{y}} L(\homo{y},X,\homo{f})$
to be maximized. Gathering the
terms depending on $\homo{y}$, the semidefinite problem dual to
(\ref{primal}) consists in finding $X$ and $\homo{f}$ maximizing
the linear function $-\homo{f}(\homo{x}^*)$ subject to
\begin{equation}\label{dual}
\begin{array}{l}
\ip{A_{00}}{X_0}+\ip{B_{00}}{X_1} = \homo{f}_0 \\
\ip{A_{10}}{X_0}+\ip{B_{10}}{X_1} = \homo{f}_1 \\
\ip{A_{01}}{X_0}+\ip{B_{01}}{X_1} = \homo{f}_2 \\
\ip{A_{\alpha}}{X_0}+\ip{B_{\alpha}}{X_1} = 0, \: |\alpha| > 1 \\
X_0 \succeq 0.
\end{array}
\end{equation}
Notice that feasibility of dual problem (\ref{dual})
amounts to the existence of an SOS representation for the affine
expression
$f(x) = s_0(x) + s_1(x)p(x)$
with $s_0(x) \in {\mathbb S}[x]_{2k}$ and
$s_1(x) \in {\mathbb R}[x]_{2(k-2)}$, with respective Gram matrices
$X_0$ and $X_1$, see e.g. \cite{l}.

Weak duality \cite[Section 5.8.1]{bv} informs us that the
optimal value of the dual objective function $-\homo{f}(\homo{x}^*)$ s.t.
constraints (\ref{dual}) is always greater than or equal to zero.
If this value is strictly positive, i.e. if there exists an $\homo{f}$ such that
$\homo{f}(\homo{x}^*) < 0$, then primal problem (\ref{primal}) is infeasible.
In turn, this implies that $\homo{x}^* \notin {\mathcal P}_2$.
Conversely, if $\homo{x}^* \in {\mathcal P}_2$ then for all $\homo{f}$
feasible for problem (\ref{dual}), i.e. for all $\homo{f} \in {\mathcal F}_k$,
it holds $\homo{f}(\homo{x}^*) \geq 0$. $\Box$.

\section{Exactness}\label{section_exact}

Given any element $\homo{f} \in {\mathcal F}$, define
the quartic
\[
p_f(x) = f(x)-p(x).
\]

\begin{lemma}\label{exact}
The first relaxation is exact, i.e. ${\mathcal P}$
is semidefinite representable as ${\mathcal P}_2$
(with at most 12 lifting variables)
if and only if, for all $\homo{f} \in {\mathcal F}$,
$p_f(x) \geq 0$ for all $x$.
\end{lemma}

{\bf Proof}:
Given any $\homo{f} \in {\mathcal F}$, the inequality $p_f(x) \geq 0$
implies that $p_f(x) = s_0(x) \in {\mathbb S}[x]_4$,
since bivariate quartics are non-negative if and only if they are polynomial SOS \cite{r}.
Since $\homo{f} \in {\mathbb P}^2$ by homogeneity
we can choose $s_1(x) = 1$
(without loss of generality) such that $f(x) = s_0(x) + s_1(x)p(x)$, and thus
$\homo{f} \in {\mathcal F}_2$. Therefore ${\mathcal F} = {\mathcal F}_2$
and hence ${\mathcal P} = {\mathcal P}_2$. The total number of liftings is
equal to 15 (the number of monomials of a trivariate quartic),
subject to 3 equality constraints, leaving 12 degrees of
freedom. $\Box$.

Define
\[
{\mathcal F}^* = \{\homo{f} \in {\mathbb P}^2 \: :\:
f(x) \geq 0 \:\: \forall x \in \partial{\mathcal P}\}.
\]
as the subset of $\mathcal F$ consisting only of
lines which are tangents to $\mathcal P$.

\begin{lemma}\label{tangent}
${\mathcal P}={\mathcal P}_2$ if and only if,
for all $\homo{f} \in {\mathcal F}^*$,
$p_f(x) \geq 0$ for all $x$.
\end{lemma}

{\bf Proof:}
Let $x^f \in \partial\mathcal P$ be the solution
of the convex (but possibly non-smooth) problem of 
maximizing the linear function $f^Tx$ subject to
the constraint that $x$ belongs to $\mathcal P$.
The line $f^T(x-x^f)=0$ is the tangent
to $\mathcal P$ at $x=x^f$. Let $\homo{f}_0=-f^Tx^f$, so that
$f(x)=f^T(x-x^f)=\homo{f}_0+\homo{f}_1x_1+\homo{f}_2x_2 \geq 0$
for all $x \in \mathcal P$, and the corresponding element $\homo{f}$
belongs to $\mathcal F$. Given such an element, assume that the
corresponding polynomial $p_f(x)$ is non-negative. Any other element
$\homo{f}^* \in \mathcal F$ such that $(\homo{f}^*_1,\:\homo{f}^*_2) = f$ 
has a constant value $\homo{f}^*_0$ which is larger than $\homo{f}_0$,
and hence the corresponding polynomial $p_{f^*}(x)$ is also
non-negative.
$\Box$

Checking the condition of Lemma \ref{tangent}
implies sweeping out over a real parameter (an angle),
and for each value, finding the tangent to $\mathcal P$
and checking non-negativity of the bivariate quartic $p_f(x)$
(e.g. by solving a semidefinite programming problem).
This can be computationally demanding,
and it makes sense to derive more tractable sufficient conditions ensuring
${\mathcal P}={\mathcal P}_2$ or ${\mathcal P} \neq {\mathcal P}_2$.

Let $\nabla p(x) \in {\mathbb R}[x]^2_3$ denote the gradient of $p(x)$.

\begin{lemma}\label{gradient}
If $\mathcal P$ is bounded, then ${\mathcal P}={\mathcal P}_2$ if and only if,
for all $\homo{f} \in {\mathcal F}^*$,
$p_f(x) \geq 0$ for all $x$ such
that $\nabla p(x) = f$.
\end{lemma} 

{\bf Proof}:
If $\mathcal P$ is bounded, then $p(x) \rightarrow -\infty$ 
and hence $p_f(x) \rightarrow +\infty$ when $\|x\| \rightarrow +\infty$.
The polynomial $p_f(x)$ then
achieves its minimum when its gradient vanishes, i.e. when $\nabla p(x)=f$.
$\Box$

Generically, there is a finite number (at most 4) of real points $x$
satisfying $\nabla p(x)=f$, and the sign of $p_f(x)$ should be tested
only at these points, which is a significant saving over assessing
global non-negativity of $p_f(x)$ as in Lemmas \ref{exact} or \ref{tangent}.
If $\mathcal P$ is not bounded,
we also have to check the sign of $p_f(x)$ at infinity.

\begin{lemma}\label{concave}
If $p(x)$ is concave, then ${\mathcal P}={\mathcal P}_2$.
\end{lemma}

{\bf Proof}:
If $p(x)$ is concave, then $p_f(x)=f(x)-p(x)$
is convex. This polynomial has a unique minimum $x^f$ when its gradient
$\nabla p_f(x) = f-\nabla p(x)$ vanishes,
with $f = (\homo{f}_1,\:\homo{f}_2)$. Given $f$, the point $x^f$
solution to $f=\nabla p(x)$ is the point along the boundary
$\partial\mathcal P$ at which the line $f^T(x-x^f)=0$ is tangent
to $\mathcal P$. Since $p(x)$ is concave and $\mathcal P$ has non-empty
interior, it cannot happen that $\nabla p(x)$ vanishes along $\partial\mathcal P$.
Therefore $\partial\mathcal P$ is necessarily smooth, and point $x^f$
always exists for any $f$. At this point, polynomial
$p_f(x)$ vanishes. Since this is the unique minimum and
$p_f(x)$ is convex, it follows that $p_f(x)$ is globally non-negative.
$\Box$

Testing concavity of $p(x)$ is equivalent to
testing negative semidefiniteness of its Hessian.
The Hessian is a 2-by-2 bivariate quadratic matrix, and
it is negative semidefinite if and only if its trace
is non-positive and its determinant is non-negative.
The first condition is trivial to test, whether the
second condition can be tested by semidefinite programming.

\begin{lemma}\label{smooth}
If $\partial\mathcal P$ is non-smooth, then ${\mathcal P} \neq {\mathcal P}_2$.
\end{lemma}

{\bf Proof}:
As in the proof of Lemma \ref{tangent}, given a direction $f$,
let $x^f \in \partial\mathcal P$ be the solution
of the convex problem of 
maximizing the linear function $f^Tx$ subject to
the constraint that $x$ belongs to $\mathcal P$.
Suppose that $\partial\mathcal P$ is non-smooth at $x^f$, which
implies that $p(x^f)=0$ and $\nabla p(x^f)=0$.
Consider the Taylor expansion of the quartic $p_f(x)=f(x)-p(x)$
around $x=x^f$, which reads 
$p_f(x) = p_f(x^f)+\nabla p_f(x^f)^T (x-x^f) + \cdots
= f^T(x-x^f) + \cdots$ where the dots indicate terms of
degree two or higher.
This polynomial has a non-zero first-order term, and hence
it cannot be globally non-negative. From Lemma \ref{exact},
it follows that ${\mathcal P} \neq {\mathcal P}_2$.
$\Box$

Lemma \ref{smooth} states that smoothness of $\partial\mathcal P$
is necessary for the first relaxation to be exact. However, it says nothing
about semidefinite representability of $\mathcal P$ in general.

Testing smoothness of $\partial\mathcal P$ is equivalent
to finding all singular points of $\mathcal C$ (this can be
done by solving the system of polynomial equations
$p(x)=\nabla p(x)=0$) and testing whether they lie in 
the interior of $\mathcal P$ or not.

\section{Examples}

\subsection{Egg}

\begin{figure}[h!]
\begin{center}
\includegraphics[width=12cm]{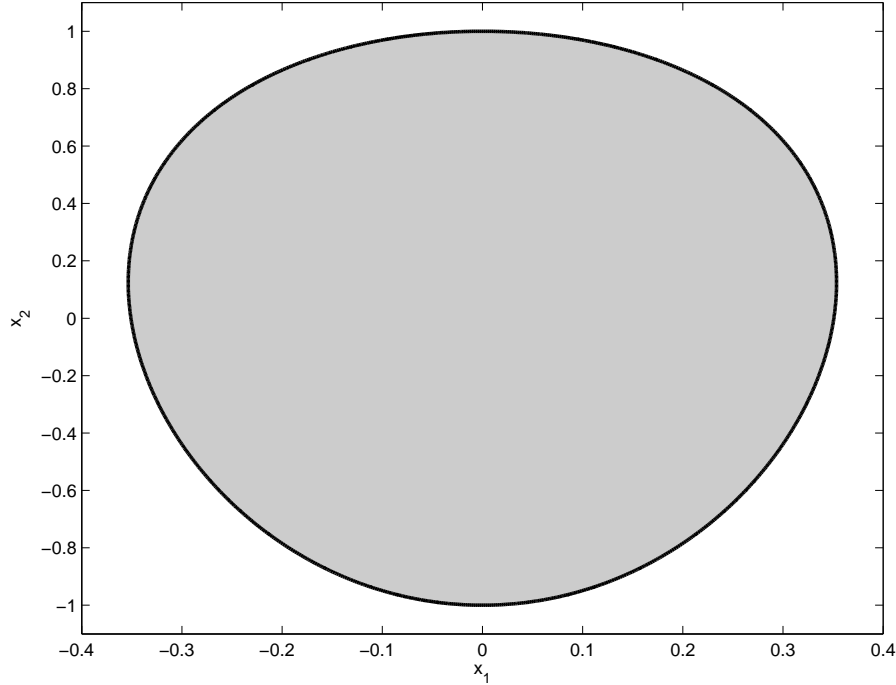}\\
\caption{Egg modeled by a 6-by-6 LMI with 12 liftings.
\label{egg_fig}}
\end{center}
\end{figure}

Let $p(x)=1-8x_1^2-(x_1^2-x_2)^2$
describe an egg curve. With the {\tt algcurves} package for Maple
we can check that curve $\mathcal C$
has genus zero with a triple singular point
$\homo{x}=(0,0,1)$ at infinity:
\begin{verbatim}
> with(algcurves):
> p:=1-8*x1^2-(x1^2-x2)^2:
> genus(p,x1,x2);
                               0
> singularities(p,x1,x2);
                     {[[0, 1, 0], 2, 3, 2]}
\end{verbatim}

The Hessian of $p(x)$ is given by
\[
\left[\begin{array}{cc}
-16-4x_2-12x_1^2 & 4x_1 \\ 4x_1 & -2
\end{array}\right]
\]
and its evaluation at $x=(0,\:-5)$ shows that it is indefinite
and hence that quartic $p(x)$ is not concave, so that Lemma \ref{concave}
cannot be applied.

Let us test the sign condition of Lemma \ref{exact}
for a given direction $\homo{f} \in {\mathcal F}^*$.
Choose e.g. the point $x^f=(0,\:1) \in {\mathcal C}$
at which $\nabla p(x^f) = (0,\:-2)$ and hence
$\homo{f} = (2,\:0,\:-2)$. We have
$p_f(x) = 2-2x_2-p(x) = 1-2x_2+8x_1^2+(x_1^2-x_2)^2$.
With YALMIP \cite{yalmip} we could find
the SOS decomposition $p_f(x)=(1-x_2+x_1^2)^2+6x_1^2$
certifying non-negativity of $p_f(x)$:
\begin{verbatim}
>> sdpvar x1 x2
>> pf=1-2*x2+8*x1^2+(x1^2-x2)^2;
>> [sol,v,Q]=solvesos(sos(pf));
>> Q{1}
ans =
    1.0000   -1.0000    1.0000         0
   -1.0000    1.0000   -1.0000         0
    1.0000   -1.0000    1.0000         0
         0         0         0    6.0000
>> eig(Q{1})'
ans =
    0.0000    0.0000    3.0000    6.0000
>> sdisplay(v{1}')
ans = 
    '1'    'x2'    'x1^2'    'x1'
>> sdisplay(clean(pf-v{1}'*Q{1}*v{1},sqrt(eps)))
0
\end{verbatim}

Sweeping over all points along $\partial{\mathcal P} = {\mathcal C}$,
we can check that the condition
of Lemma \ref{exact} is satisfied, and 
the convex set ${\mathcal P} = {\mathcal P}_2$
is represented on Figure \ref{egg_fig}.

\subsection{Bean}\label{bean}

\begin{figure}[h!]
\begin{center}
\includegraphics[width=12cm]{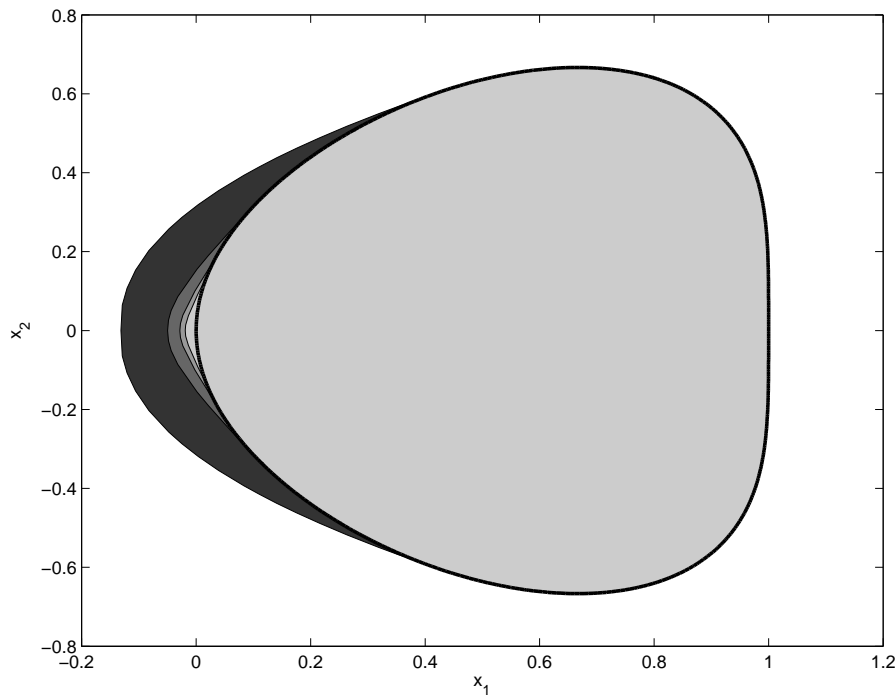}\\
\caption{Bean curve (thick line) and its first four
embedded outer semidefinite approximations (thin lines and
shaded regions).\label{bean_fig}}
\end{center}
\end{figure}

Let $p(x)=x_1(x_1^2+x_2^2)-x_1^4-x_1^2x_2^2-x_2^4$
define the so-called bean quartic. The curve has genus zero,
and a triple singular point at the origin $\homo{x}=(1,0,0)$.

At this point the gradient $\nabla p(x)$ vanishes, whereas
the curve $\mathcal P$ has a tangent $f(x)=x_1$.
Polynomial $p_f(x)=x_1-p(x)$ has a non-zero linear term,
and hence it cannot be non-negative. By Lemma \ref{smooth}
we have ${\mathcal P}$ strictly included in ${\mathcal P}_2$.

Inspection of the constant and linear terms in the expression (\ref{fk})
arising in the definition of ${\mathcal P}_k$ in Lemma \ref{lemma_dual}
shows that actually $\mathcal P$ is strictly included
in ${\mathcal P}_k$ for all $k$.

On Figure \ref{bean_fig} we see embedded semidefinite representable
sets ${\mathcal P}_k$ for $k=2,3,4,5$ (thin lines
and shaded regions) and the convex set ${\mathcal P}$ (thick line).
It seems that $\mathcal P$ is smooth but actually there is
a singularity at the origin. We see the global consequences
of the pointwise singularity on the shape of the
sets ${\mathcal P}_k$.

For interested readers, semidefinite sets ${\mathcal P}_k$ can
be visualized with the following Matlab script, mixing features
from GloptiPoly 3 \cite{g3} and YALMIP \cite{yalmip}:
\begin{verbatim}
mpol x 2
p = x(1)*(x(1)^2+x(2)^2)-x(1)^4-x(1)^2*x(2)^2-x(2)^4;
k = 3; % relaxation order
P = msdp(p==0, k);
[F,h,y] = myalmip(P);
plot(F,y(1:2)); % projection on first degree moments 
\end{verbatim}

The impact in polynomial optimization
of the non-exactness of semidefinite relaxations
can be observed with the help of the following Matlab script
using GloptiPoly 3:
\begin{verbatim}
bounds = [];
for k = 2:10
 P = msdp(min(x(1)), p==0, k);
 [status,obj] = msol(P);
 bounds = [bounds obj];
end
\end{verbatim}
We obtain the following sequence of lower bounds on
the minimum value of $x_1$ such that $p(x)$ vanishes:
$-0.1315$, $-0.02915$, $-0.009705$, $-0.01022$,
$-0.009319$, $-0.009320$, $-0.009646$, $-0.009346$, $-0.009371,\ldots$
We expect the sequence to converge from below to zero, the genuine minimum,
but numerically, it stagnates around $-9\cdot 10^{-3}$.

\subsection{Water drop}\label{waterdrop}

\begin{figure}[h!]
\begin{center}
\includegraphics[width=12cm]{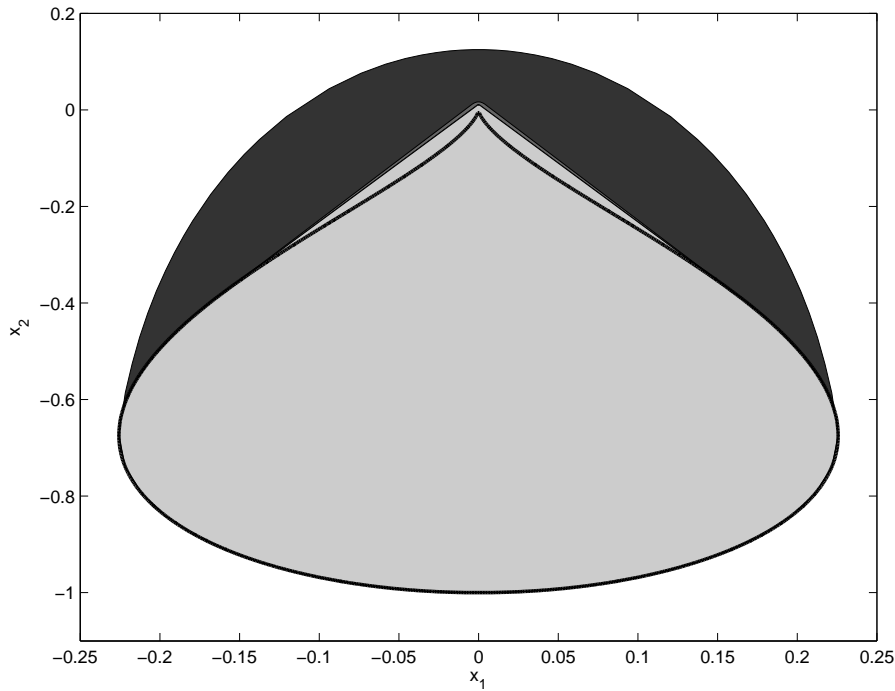}\\
\caption{Water drop quartic (thick line)
and its first four embedded semidefinite approximations
(thin lines and shaded regions).\label{waterdrop_fig}}
\end{center}
\end{figure}

Let $p(x)=-x_1^2-x_2^3-(x_1^2+x_2^2)^2$ define
a water drop quartic. The curve has genus two, with
a singular point (a cusp) at the origin.
On Figure \ref{waterdrop_fig} we see semidefinite representable
embedded sets ${\mathcal P}_k$ for $k=2,3,4,5$ (thin lines
and shaded regions) and the non convex curve ${\mathcal C}$ (thick line).
The set ${\mathcal P}_5$ and the convex hull ${\mathcal P} =
\mathrm{conv}\:{\mathcal C}$ are almost undistinguishable.
However, as for Example \ref{bean}, it can be shown that
${\mathcal P}_k$ cannot be equal to $\mathcal P$
for $k$ finite.

\subsection{Lemniscate}\label{lemniscate}

\begin{figure}[h!]
\begin{center}
\includegraphics[width=12cm]{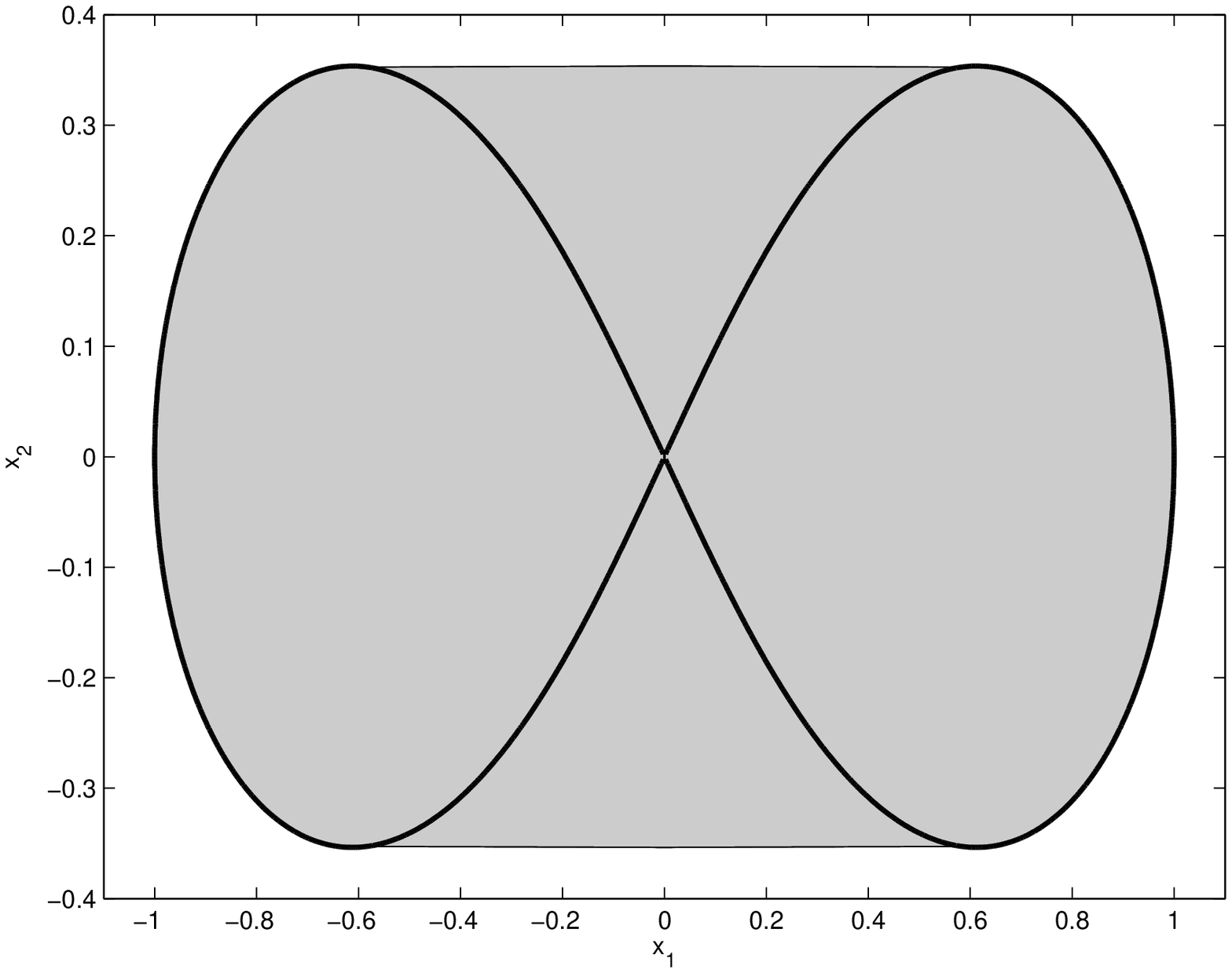}\\
\caption{Lemniscate (thick line) and its convex hull
modeled by a 6-by-6 LMI with 12 liftings.
\label{lemniscate_fig}}
\end{center}
\end{figure}

Let $p(x)=x_1^2-x_2^2-(x_1^2+x_2^2)^2$ define a lemniscate,
a curve of genus zero, with singular points at the origin
and at the infinite complex points $\homo{x}=(0,\pm i,1)$.
Even though $\mathcal C$ is singular, the boundary
$\partial\mathcal P$ is smooth.

Sweeping over all directions $f$ indicates that $p_f(x)$
is always non-negative, and hence that ${\mathcal P}={\mathcal P}_2$,
see Figure \ref{lemniscate}. Here, the singularity
of $\mathcal C$ is in the interior of $\mathcal P$,
and it does not prevent the first relaxation
to be exact.

\subsection{Folium}\label{folium}

\begin{figure}[h!]
\begin{center}
\includegraphics[width=12cm]{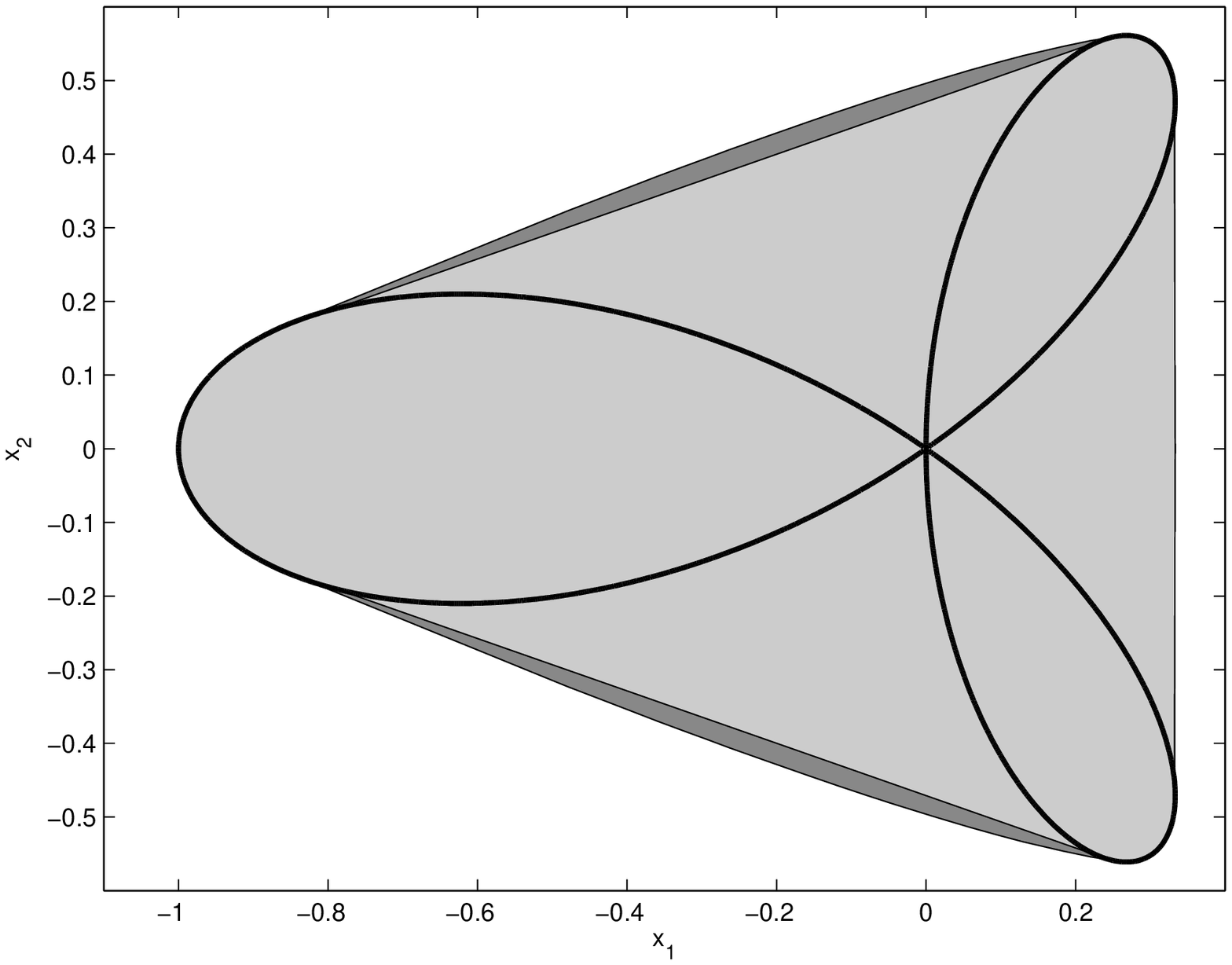}\\
\caption{Folium curve
and its first two embedded outer semidefinite approximations
(thin lines and shaded regions).\label{folium_fig}}
\end{center}
\end{figure}

Let $p(x)=-x_1(x_1^2-2x_2^2)-(x_1^2+x_2^2)^2$ define a folium,
a curve of genus zero, with a triple singular point at the origin.
As in Example \ref{lemniscate}, the singularity of $\mathcal C$
is the interior of $\mathcal P$, hence the boundary $\partial\mathcal P$
is smooth, and we may expect that the first relaxation is exact.

However, along the direction $\homo{f}=(1,\:3/4,\:3\sqrt{2}/2)$
corresponding to a bitangent of $\mathcal C$ (obtained computationally
by finding singular points of the dual curve),
polynomial $p_f(x)=f(x)-p(x)$ has a strictly negative minimum
achieved at $x \approx (0.2040,\:-0.8762)$.

On Figure \ref{folium_fig} we see curve $\mathcal C$ (thick line)
and semidefinite sets ${\mathcal P}_2$ (exterior thin line, dark shaded
region) and ${\mathcal P}_3$ (interior thin line, shaded region).
We observe that ${\mathcal P} = \mathrm{conv}\:{\mathcal C}$
is strictly included in ${\mathcal P}_2$,
whereas, apparently, ${\mathcal P}_3={\mathcal P}$ (but
we are not able to prove this identity).

\subsection{Smooth and convex}

\begin{figure}[h!]
\begin{center}
\includegraphics[width=12cm]{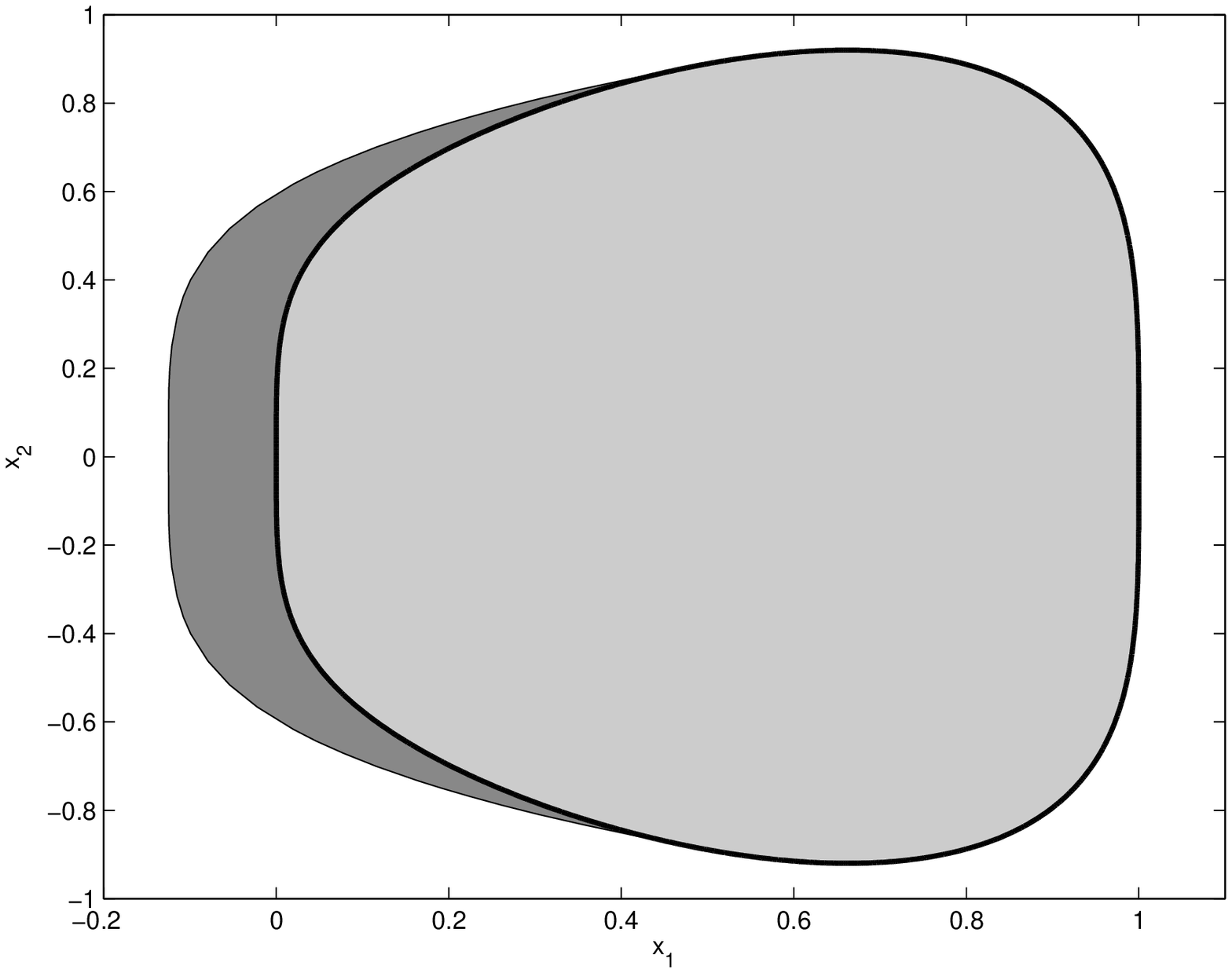}\\
\caption{Smooth convex quartic
and its first two embedded outer semidefinite approximations
(thin lines and shaded regions).\label{smoothconvex_fig}}
\end{center}
\end{figure}

Let $p(x)=x_1+x_1^2-2x_1^4-x_2^4$ define a convex
quartic curve of genus three. Note that $p(x)$ is not
concave. There is no singularity, so we may expect that
the first relaxation is exact. However, for the
tangent $f(x)=x_1$, polynomial $p_f(x)=f(x)-p(x)$
is not non-negative (consider e.g. its sign along
the line $x_2=0$) and hence ${\mathcal P} \neq {\mathcal P}_2$.
As in Example \ref{folium}, apparently
${\mathcal P}={\mathcal P}_3$.

This example was prepared with the help of Bill Helton
and Jiawang Nie.

\subsection{Fermat}\label{tv}

The Fermat quartic $p(x)=1-x_1^4-x_2^4$ (also called
the TV-screen quartic) is smooth,
so its convex hull is semidefinite representable as
a 6-by-6 LMI with 12 liftings. However, inspection reveals
that there is a semidefinite representation with only 2 liftings:
\[
{\mathcal P} = \{\homo{x} \in {\mathbb P}^2 \: :\: \exists \homo{y} \in 
{\mathbb P} \: :\: \left[\begin{array}{cccccc}
\homo{x}_0+\homo{y}_0 & \homo{y}_1 \\
\homo{y}_1 & \homo{x}_0-\homo{y}_0 \\
& & \homo{x}_0 & \homo{x}_1 \\
& & \homo{x}_1 & \homo{y}_0 \\
& & & & \homo{x}_0 & \homo{x}_2 \\
& & & & \homo{x}_2 & \homo{y}_1
\end{array}\right] \succeq 0\}.
\]
However, we do not know how to derive this reduced
representation with 2 liftings from the generic
representation with 12 liftings.
 
\section{Rational curves}\label{rat}

In this section we restrict
the class of $\mathcal C$ to algebraic curves of genus zero
\cite{h},
i.e. curves which admit a polynomial parametrization
\begin{equation}\label{param}
{\mathcal C} =
\{\homo{x} \in {\mathbb P}^2 \: :\: \exists\: \homo{t} \in {\mathbb P} \: :\: 
\homo{x}_i = \homo{p}_i(\homo{t}), \: i=0,1,2\}
\end{equation}
with $\homo{p}_i(\homo{t})$ bivariate quartic forms
in $\homo{t} \in {\mathbb P}$.
Given $p(x)$ in implicit representation (\ref{curve}),
there are algorithms to compute the $\homo{p}_i(\homo{t})$ in the above
explicit representation, see e.g. the Maple package {\tt algcurves}
for an implementation.

Let us define the Hankel moment matrix
\[
M_2(\homo{y})
= \sum_{\alpha}^4 \homo{y}_{\alpha} H_{\alpha} = \left[\begin{array}{cccc}
\homo{y}_0 & \homo{y}_1 & \homo{y}_2 \\
\homo{y}_1 & \homo{y}_2 & \homo{y}_3 \\
\homo{y}_2 & \homo{y}_3 & \homo{y}_4
\end{array}\right]
\]
and, given the quartic $\homo{p}(\homo{t})=\sum_{\alpha}^4
\homo{p}_{\alpha} \homo{t}^{\alpha}$,
the localising linear form
\[
M_0(\homo{p}\homo{y}) = \sum_{\alpha}^4 \homo{p}_{\alpha} \homo{y}_{\alpha}.
\]

\begin{lemma}\label{exactrat}
${\mathcal P} = \{\homo{x} \in {\mathbb P}^2 \: :\: \exists\: \homo{y} \in {\mathbb P}^4 \: :\:
\homo{x}_i = M_0(\homo{p}_i \homo{y}), \: i=0,1,2,\:
M_2(\homo{y}) \succeq 0\}.$
\end{lemma}

{\bf Proof:}
We closely follow the proof of Lemma \ref{exact}.
Given $\homo{f} \in {\mathcal F}$ and $\homo{x} \in {\mathcal C}$,
consider the bivariate form $\homo{g}(\homo{t})
= \homo{f}^T\homo{x} = \sum_{i=0}^2 \homo{f}_i \homo{x}_i =
\sum_{i=0}^2 \homo{f}_i \homo{p}_i(\homo{t}) =
\sum_{\alpha}^4 (\homo{g}^T_{\alpha}\homo{f}) \homo{t}^{\alpha}$
with $\homo{t} \in {\mathbb P}$ and coefficients $\homo{g}_{\alpha} \in {\mathbb R}^3$.
This non-negative bivariate form is always SOS
\cite{r}, hence there exists a 3-by-3 matrix $X$
such that
\[
\begin{array}{l}
\mathrm{trace}(H_{\alpha} X) = \homo{g}^T_{\alpha}\homo{f}, \: \alpha=0,1,\ldots,4 \\
X \succeq 0,
\end{array}
\] 
hence the dual formulation of Lemma \ref{exactrat}.$\Box$

An alternative, more direct proof suggested by Roland Hildebrand,
consists in viewing $\mathcal P$ as the image through a linear mapping
of the convex hull of a Veronese variety, namely the image of the nonlinear map
sending $\homo{t} \in {\mathbb P}$ into $[\homo{t}_0^4 \:\: \homo{t}_0^3\homo{t}_1 \:\: \homo{t}_0^2\homo{t}_1^2 \:\:
\homo{t}_0\homo{t}_1^3 \:\: \homo{t}_1^4]
\in {\mathbb P}^4$, see \cite{h}. This Veronese variety is also
called sometimes the moment curve, and its convex hull is indeed
the cone of positive semidefinite Hankel matrices.

Jean-Bernard Lasserre informed me that Pablo Parrilo presented a related
semidefinite representation of the convex hull of rational plane curves
at a workshop at Banff, Canada, in October 2006. At the time of
writing of these notes (March 2008), the result is not available
in electronic or printed form, however.

Note that the relations $\homo{x}_i = M_0(\homo{p}_i\homo{y})$ in Lemma \ref{exactrat}
form a consistent linear system of 3 equations with 5 indeterminates,
so the number of lifting variables
can always be reduced to 5-3=2.

\subsection{Folium revisited}\label{folium2}

Consider again the folium quartic of Example \ref{folium}.
With the {\tt algcurves} package of Maple, we obtain
the following rational parametrization:
$p_0(t)=1+2t_1^2+t_1^4$, $p_1(t)=-1+2t_1^2$,
$p_2(t)=-t_1+2t_1^3$.

The lifting variables in the representation of Lemma \ref{exactrat}
satisfy the linear system of equations
$\homo{x}_0 = \homo{y}_0 + 2\homo{y}_2 + \homo{y}_4$, $\homo{x}_1 = -\homo{y}_0 + 2\homo{y}_2$,
$\homo{x}_2 = -\homo{y}_1 + 2\homo{y}_3$.
From this we derive $\homo{y}_2 = \frac{1}{2}(\homo{x}_1+\homo{y}_0)$, 
$\homo{y}_3 = \frac{1}{2}(\homo{x}_2+\homo{y}_1)$, $\homo{y}_4 = \homo{x}_0-\homo{x}_1-2\homo{y}_0$
that we can report in the Hankel matrix constraint
to produce a semidefinite representation
of the convex hull of $\mathcal C$ with 2 liftings:
\[
{\mathcal P} = \{\homo{x} \in {\mathbb P}^2 \: :\: \exists \homo{y} \in {\mathbb P} \: :\:
\left[\begin{array}{ccc}
2\homo{y}_0 & 2\homo{y}_1 & \homo{x}_1+\homo{y}_0 \\
2\homo{y}_1 & \homo{x}_1+\homo{y}_0 & \homo{x}_2+\homo{y}_1 \\
\homo{x}_1+\homo{y}_0 & \homo{x}_2+\homo{y}_1 & 2\homo{x}_0-2\homo{x}_1-4\homo{y}_0 
\end{array}\right] \succeq 0\}.
\]

\subsection{Bean revisited}

Consider again the bean quartic of Example \ref{bean}.
A rational parametrization (\ref{param})
is given by $p_0(t)=1+t_1^2+t_1^4$, $p_1(t)=1+t_1^2$
and $p_2(t)=t_1+t_1^3$, from which follows 
a semidefinite representation with 2 liftings:
\[
{\mathcal P} = \{\homo{x} \in {\mathbb P}^2 \: :\: \exists \homo{y} \in {\mathbb P} \: :\:
\left[\begin{array}{ccc}
\homo{y}_0 & \homo{y}_1 & \homo{x}_1-\homo{y}_0 \\
\homo{y}_1 & \homo{x}_1-\homo{y}_0 & \homo{x}_2-\homo{y}_1 \\
\homo{x}_1-\homo{y}_0 & \homo{x}_2-\homo{y}_1 & \homo{y}_0-\homo{y}_1
\end{array}\right] \succeq 0\}.
\]

\section{Conclusion}

This note investigated semidefinite representations
of convex plane quartics, and more specifically the exactness
of the first semidefinite relaxation in Lasserre's hierarchy.
Also described was an elementary
exact semidefinite representation of the convex hull
of rationally parametrized quartics.
Exactness conditions followed from the well-known fact that non-negative
polynomials can be represented as sum-of-squares in the bivariate
quartic case and in the univariate case. 
It follows that the exactness
result of Section \ref{rat} is valid for rationally parametrizable
curves of arbitrary degree and dimension, but the exactness
result of Section \ref{section_exact} is limited
to plane quartics. Also unclear is what kind of conditions
should be enforced to ensure exactness of the second, third,
and in general higher-order relaxations.

In \cite{l06}, Lasserre proposed sufficient conditions for
convex semialgebraic sets to be semidefinite representable.
The conditions are algebraic in nature, strongly connected
with the polynomials used to model the set, and related
with the particular SOS representation of Lemma \ref{lemma_dual}.
We also notice that Lemma \ref{concave} can be found in Example 3
in \cite{l06} where it is proved using Karush-Kuhn-Tucker
optimality conditions.
In \cite{hn}, Helton and Nie derived sufficient conditions
in terms of negative definiteness of the Hessian along the tangent
space along the boundary of the set, provided the gradient does not vanish
along this boundary. In contrast with these general statements,
our focus in this note was more on computational aspects
and explicit examples, the driving force being that if we do not understand
well the simplest non-trivial case (plane quartics) it is likely
that we will not understand more complicated configurations.
It is expected that our examples can provide further motivation
for studying alternative semidefinite representations, see e.g.
\cite{nie}. For instance, we are not aware of any exact semidefinite
representation of the convex hull of the water drop quartic
of Example \ref{waterdrop}.

As illustrated in Example \ref{tv}, a given quartic may have 
different semidefinite representations with a different number
of lifting variables. Given a
representation, it could be interesting to design a systematic
algorithm to remove redundant lifting variables. Similarly, the
problem of finding a representation with a minimum number
of lifting variables seems to be open.

Since projections of LMIs are convex semialgebraic
sets, and the essential difficulty when building semidefinite
representations seems to be singularities (points at which the gradient
vanishes), one may be tempted to conjecture that convex regions delimited
by higher degree smooth curves admit an exact semidefinite representation. 
Helton and Nie go even farther in the conclusion of \cite{hn}
by conjecturing that every convex semialgebraic set is
semidefinite representable.

\section*{Appendix}

In this paper we use projective
spaces ${\mathbb P}^k$ over the field $\mathbb R$, together with
affine spaces ${\mathbb R}^k$. By projective space
${\mathbb P}^k$, we mean the set of all one-dimensional
subspaces of ${\mathbb R}^{k+1}$. Equivalently, ${\mathbb P}^k$
is the quotient space of equivalence classes of ${\mathbb R}^{k+1}-{0}$
under the equivalence relation $(\homo{x}_0,\homo{x}_1,\ldots,\homo{x}_k)
\sim (a\homo{x}_0,a\homo{x}_1,\ldots,a\homo{x}_k)$ for all nonzero $a \in {\mathbb R}$.
Projective space ${\mathbb P}^k$ is a compact space under
the Zariski topology where a closed set is defined as the
zero set of homogeneous polynomials. When $\homo{x}_0\neq 0$, it holds 
$(\homo{x}_0,\homo{x}_1,\ldots,\homo{x}_k) \sim (1,\homo{x}_1/\homo{x}_0,\ldots,\homo{x}_k/\homo{x}_0)$, and we can identify
a point $(x_1,\ldots,x_k) \in {\mathbb R}^k$ with
a point $(1,x_1,\ldots,x_k) \in {\mathbb P}^k$.
Then the affine space ${\mathbb R}^k$ is the open subset of
the projective space ${\mathbb P}^k$
defined by $\homo{x}_0 \neq 0$. Points $(0,\homo{x}_1,\ldots,\homo{x}_k)$ with
$\homo{x}_0=0$ corresponds to points at infinity, and
${\mathbb P}^k$ can be also viewed as the affine space ${\mathbb R}^k$
extended with points at infinity. See e.g. \cite[Chapter 1]{h}
for an elementary introduction.

\section*{Acknowledgments}

This work benefited from exchanges with Roland Hildebrand, Bill Helton,
Jean-Bernard Lasserre and Jiawang Nie.

\end{document}